# NEAR-IDEAL MODEL SELECTION BY $\ell_1$ MINIMIZATION

By Emmanuel J. Candès[1] and Yaniv Plan

*California Institute of Technology*

We consider the fundamental problem of estimating the mean of a vector $y = X\beta + z$, where $X$ is an $n \times p$ design matrix in which one can have far more variables than observations, and $z$ is a stochastic error term—the so-called "$p > n$" setup. When $\beta$ is sparse, or, more generally, when there is a sparse subset of covariates providing a close approximation to the unknown mean vector, we ask whether or not it is possible to accurately estimate $X\beta$ using a computationally tractable algorithm.

We show that, in a surprisingly wide range of situations, the lasso happens to nearly select the best subset of variables. Quantitatively speaking, we prove that solving a simple quadratic program achieves a squared error within a logarithmic factor of the ideal mean squared error that one would achieve with an *oracle* supplying perfect information about which variables should and should not be included in the model. Interestingly, our results describe the average performance of the lasso; that is, the performance one can expect in an vast majority of cases where $X\beta$ is a sparse or nearly sparse superposition of variables, but not in all cases.

Our results are nonasymptotic and widely applicable, since they simply require that pairs of predictor variables are not too collinear.

## 1. Introduction.

One of the most common problems in statistics is to estimate a mean response $X\beta$ from the data $y = (y_1, y_2, \ldots, y_n)$ and the linear model

$$(1.1) \qquad\qquad y = X\beta + z,$$

where $X$ is an $n \times p$ matrix of explanatory variables, $\beta$ is a $p$-dimensional parameter of interest and $z = (z_1, \ldots, z_n)$ is a vector of independent stochastic

---

Received June 2008.

[1]Supported in part by a National Science Foundation Grant CCF-515362 and by the 2006 Waterman Award (NSF).

errors. Unless specified otherwise, we will assume that the errors are Gaussian with $z_i \sim \mathcal{N}(0, \sigma^2)$, but this is not really essential, as our results and methods can easily accommodate other types of distribution. We measure the performance of any estimator $X\hat{\beta}$ with the usual squared Euclidean distance $\|X\beta - X\hat{\beta}\|_{\ell_2}^2$, or with the mean-squared error, which is simply the expected value of this quantity.

In this paper, although this is not a restriction, we are primarily interested in situations where there are as many explanatory variables as observations, or more—the so-called, and now widely popular, "$p > n$" setup. In such circumstances, however, it is often the case that a relatively small number of variables have substantial explanatory power, so that, to achieve accurate estimation, one needs to select the "right" variables and determine which components $\beta_i$ are not equal to zero. A standard approach is to find $\hat{\beta}$ by solving

$$(1.2) \qquad \min_{b \in \mathbb{R}^p} \frac{1}{2}\|y - Xb\|_{\ell_2}^2 + \lambda_0 \sigma^2 \|b\|_{\ell_0},$$

where $\|b\|_{\ell_0}$ is the number of nonzero components in $b$. In other words, the estimator (1.2) achieves the best trade-off between the goodness of fit and the complexity—in this case, the number of variables included—of the model. Popular selection procedures such as AIC, $C_p$, BIC and RIC are all of this form, with different values of the parameter ($\lambda_0 = 1$ in AIC [1, 19], $\lambda_0 = \frac{1}{2}\log n$ in BIC [24], and $\lambda_0 = \log p$ in RIC [14]). It is known that these methods perform well both empirically and theoretically (see [14] and [2, 4] and the many references therein). Having said this, the problem, of course, is that these "canonical selection procedures" are highly impractical. Solving (1.2) is, in general, NP-hard [22] and, to the best of our knowledge, requires exhaustive searches over all subsets of columns of $X$, a procedure which is clearly combinatorial in nature and has exponential complexity, since, for $p$ of size about $n$, there are about $2^p$ such subsets.

In recent years, several methods based on $\ell_1$ minimization have been proposed to overcome this problem. The most well-known is probably the lasso [26], which replaces the nonconvex $\ell_0$ norm in (1.2) with the convex $\ell_1$ norm $\|b\|_{\ell_1} = \sum_{i=1}^p |b_i|$. The lasso estimate $\hat{\beta}$ is defined as the solution to

$$(1.3) \qquad \min_{b \in \mathbb{R}^p} \frac{1}{2}\|y - Xb\|_{\ell_2}^2 + \lambda\sigma\|b\|_{\ell_1},$$

where $\lambda$ is a regularization parameter that essentially controls the sparsity (or the complexity) of the estimated coefficients (see [23] and [11] for exactly the same proposal). In contrast to (1.2), the optimization problem (1.3) is a quadratic program that can be solved efficiently. It is known that the lasso performs well in some circumstances. Further, there is also an emerging literature on its theoretical properties [3, 5, 6, 15, 16, 20, 21, 28, 29, 30] showing that, in some special cases, the lasso is effective.

In this paper, we will show that the lasso works provably well in a surprisingly broad range of situations. We establish that, under minimal assumptions guaranteeing that the predictor variables are not highly correlated, the lasso achieves a squared error nearly as good as if one had an oracle supplying perfect information about which $\beta_i$'s were nonzero. Continuing in this direction, we also establish that the lasso correctly identifies the true model with very large probability, provided that the amplitudes of the nonzero $\beta_i$ are sufficiently large.

1.1. *The coherence property.* Throughout the paper, we will assume that, without loss of generality, the matrix $X$ has unit-normed columns, as one can otherwise always rescale the columns. We denote, by $X_i$, the $i$th column of $X$ ($\|X_i\|_{\ell_2} = 1$) and introduce the notion of coherence, which essentially measures the maximum correlation between unit-normed predictor variables and is defined by

$$(1.4) \qquad \mu(X) = \sup_{1 \leq i < j \leq p} |\langle X_i, X_j \rangle|.$$

In other words, the coherence is the maximum inner product between any two distinct columns of $X$. It follows that, if the columns have zero mean, the coherence is just the maximum correlation between pairs of predictor variables.

We will be interested in problems in which the variables are not highly collinear or redundant.

DEFINITION 1.1 (*Coherence property*). A matrix $X$ is said to obey the coherence property if

$$(1.5) \qquad \mu(X) \leq A_0 \cdot (\log p)^{-1},$$

where $A_0$ is some positive numerical constant.

A matrix obeying the coherence property is a matrix in which the predictors are not highly collinear. This is a mild assumption. Suppose that $X$ is a Gaussian matrix with i.i.d. entries whose columns are subsequently normalized. The coherence of $X$ is about $\sqrt{(2 \log p)/n}$, so that such matrices trivially obey the coherence property, unless $n$ is ridiculously small [i.e., of the order of $(\log p)^3$]. We will give other examples of matrices obeying this property later in the paper, and we will soon contrast this assumption with what is traditionally assumed in the literature.

1.2. *Sparse model selection.* We begin by discussing the intuitive case, where the vector $\beta$ is sparse, before extending our results to a completely general case. The basic question we would like to address here is, how well can one estimate the response $X\beta$, when $\beta$ happens to have only $S$ nonzero components? From now on, we call such vectors $S$-*sparse*.

First and foremost, we would like to emphasize that, in this paper, we are interested in quantifying the performance one can expect from the lasso in an overwhelming majority of cases. This viewpoint needs to be contrasted with an analysis concentrating on the worst case performance; when the focus is on the worst case scenario, one would study very particular values of the parameter $\beta$ for which the lasso does not work well. This is not our objective; as an aside, this will enable us to show that one can reliably estimate the mean response $X\beta$ under much weaker conditions than what is currently known.

Our point of view emphasizes the average performance (or the performance one could expect in a large majority of cases); thus, we need a statistical description of sparse models. To this end, we introduce the *generic $S$-sparse model*, which is defined as follows:

1. The support $I \subset \{1, \ldots, p\}$ of the $S$ nonzero coefficients of $\beta$ is selected uniformly at random.
2. Conditional on $I$, the signs of the nonzero entries of $\beta$ are independent and equally likely to be $-1$ or $1$.

We make no assumption on the amplitudes. In some sense, this is the simplest statistical model one could think of; it says, simply, that all subsets of a given cardinality are equally likely, and that the signs of the coefficients are equally likely. In other words, one is not biased towards certain variables, nor do we have any reason to believe a priori that a given coefficient is positive or negative.

Our first result is that for most $S$-sparse vectors $\beta$, the lasso is provably accurate. Throughout, $\|X\|$ refers to the operator norm of the matrix $A$ (the largest singular value).

THEOREM 1.2. *Suppose that $X$ obeys the coherence property, and assume that $\beta$ is taken from the generic $S$-sparse model. Suppose that $S \leq c_0 p/[\|X\|^2 \log p]$ for some positive numerical constant $c_0$. Then, the lasso estimate (1.3) computed with $\lambda = 2\sqrt{2 \log p}$ obeys*

$$
(1.6) \qquad \|X\beta - X\hat{\beta}\|_{\ell_2}^2 \leq C_0 \cdot (2 \log p) \cdot S \cdot \sigma^2
$$

*with probability at least $1 - 6p^{-2\log 2} - p^{-1}(2\pi \log p)^{-1/2}$. The constant $C_0$ may be taken as $8(1 + \sqrt{2})^2$.*

For simplicity, we have chosen $\lambda = 2\sqrt{2 \log p}$, but one could take any $\lambda$ of the form $\lambda = (1 + a)\sqrt{2 \log p}$ with $a > 0$. Our proof indicates that, as $a$ decreases, the probability with which (1.6) holds decreases, but the constant $C_0$ also decreases. Conversely, as $a$ increases, the probability with which (1.6) holds increases, but the constant $C_0$ also increases.

Theorem 1.2 asserts that one can estimate $X\beta$ with nearly the same accuracy as if one knew ahead of time which $\beta_i$'s were nonzero. To see why this is true, suppose that the support $I$ of the true $\beta$ was known. In this ideal situation, we would presumably estimate $\beta$ by regressing $y$ onto the columns of $X$, with indices in $I$, and construct

$$(1.7) \quad \beta^{\star} = \arg\min_{b \in \mathbb{R}^p} \|y - Xb\|_{\ell^2}^2 \quad \text{subject to} \quad b_i = 0 \qquad \text{for all } i \notin I.$$

It is a simple calculation to show that this ideal estimator (it is ideal, because we would not know the set of nonzero coordinates) achieves[2]

$$(1.8) \qquad\qquad \mathbb{E}\|X\beta - X\beta^{\star}\|_{\ell_2}^2 = S \cdot \sigma^2.$$

Hence, one can see that (1.6) is optimal up to a factor proportional to $\log p$. It is also known that one cannot, in general, hope for a better result; the log factor is the price we need to pay for not knowing ahead of time which of the predictors are actually included in the model.

The assumptions of our theorem are pretty mild. Roughly speaking, if the predictors are not too collinear, and if $S$ is not too large, then the lasso works most of the time. An important point here is that the restriction on the sparsity can be very mild. We give two examples to illustrate our purpose:

- *Random design.* Imagine, as before, that the entries of $X$ are i.i.d. $\mathcal{N}(0, 1)$ and then normalized. Then, the operator norm of $X$ is sharply concentrated around $\sqrt{p/n}$, so that our assumption essentially reads $S \leq c_0 n / \log p$. Expressed in a different way, $\beta$ does not have to be sparse at all. It has to be smaller, of course, than the number of observations, but not by a very large margin.

  Similar conclusions would apply to many other types of random matrices.

- *Signal estimation.* A problem that has attracted quite a bit of attention in the signal processing community is that of recovering a signal that has a sparse expansion as a superposition of spikes and sinusoids. Here, we have noisy data $y$

$$(1.9) \qquad\qquad y(t) = f(t) + z(t), \qquad t = 1, \ldots, n,$$

---

[2]One could also develop a similar estimate with high probability, but we find it simpler and more elegant to derive the performance in terms of expectation.

about a digital signal $f$ of interest, which is expressed as the "time-frequency" superposition

$$(1.10) \qquad f(t) = \sum_{k=1}^{n} \alpha_k^{(0)} \delta(t-k) + \sum_{k=1}^{n} \alpha_k^{(1)} \varphi_k(t);$$

$\delta$ is a Dirac or spike obeying $\delta(t) = 1$ if $t = 0$, and 0 otherwise. $(\varphi_k(t))_{1 \le k \le n}$ is an orthonormal basis of sinusoids. The problem (1.9) is of the general form (1.1) with $X = [I_n F_n]$ in which $I_n$ is the identity matrix, $F_n$ is the basis of sinusoids (a discrete cosine transform) and $\beta$ is the concatenation of $\alpha^{(0)}$ and $\alpha^{(1)}$. Here, $p = 2n$ and $\|X\| = \sqrt{2}$. Also, $X$ obeys the coherence property if $n$ or $p$ is not too small, since $\mu(X) = \sqrt{2/n} = 2/\sqrt{p}$.

Hence, if the signal has a sparse expansion with fewer than on the order of $n/\log n$ coefficients, then the lasso achieves a quality of reconstruction that is essentially as good as what could be achieved if we knew in advance the precise location of the spikes and the exact frequencies of the sinusoids.

This fact extends to other pairs of orthobases and to general overcomplete expansions, as we will explain later.

In our two examples, the condition of Theorem 1.2 is satisfied for $S$ as large as on the order of $n/\log p$; that is, $\beta$ may have a large number of nonzero components. The novelty here is that the assumptions on the sparsity level $S$, and on the correlation between predictors, are very realistic. This is different from the available literature, which typically requires a much lower bound on the coherence or a much lower sparsity level (see Section 4 for a comprehensive discussion). In addition, many published results assume that the entries of the design matrix $X$ are sampled from a probability distribution (e.g., are i.i.d. samples from the standard normal distribution), which we are not assuming here (one could of course specialize our results to random designs as discussed above). Hence, we do not simply prove that in some idealized setting the lasso would do well, but that it has a very concrete edge in practical situations—as shown empirically in a great number of works.

An interesting fact is that one cannot expect (1.6) to hold for all models, as one can construct simple examples of incoherent matrices and special $\beta$ for which the lasso does not select a good model (see Section 2). In this sense, (1.6) can be achieved on the average—or better, in an overwhelming majority of cases—but not in all cases.

### 1.3. *Exact model recovery.*

Suppose, now, that we are interested in estimating the set $I = \{i : \beta_i \ne 0\}$. Then, we show that, if the values of the nonvanishing $\beta_i$'s are not too small, then the lasso correctly identifies the "right" model.

THEOREM 1.3. *Let $I$ be the support of $\beta$, and suppose that*

$$\min_{i \in I} |\beta_i| > 8\sigma\sqrt{2\log p}.$$

*Then, under the assumptions of Theorem 1.2, the lasso estimate with $\lambda = 2\sqrt{2\log p}$ obeys*

(1.11) $$\operatorname{supp}(\hat{\beta}) = \operatorname{supp}(\beta) \quad and$$

(1.12) $$\operatorname{sgn}(\hat{\beta}_i) = \operatorname{sgn}(\beta_i) \quad for\ all\ i \in I,$$

*with probability at least $1 - 2p^{-1}((2\pi\log p)^{-1/2} + |I|p^{-1}) - O(p^{-2\log 2})$.*

In other words, if the nonzero coefficients are significant in the sense that they stand above the noise, then the lasso identifies all the variables of interest and only these. Further, the lasso correctly estimates the signs of the corresponding coefficients. Again, this does not hold for all $\beta$'s, as shown in the example of Section 2, but for a wide majority.

Our condition says that the amplitudes must be larger than a constant times the noise level times $\sqrt{2\log p}$, which is sharp, modulo a small multiplicative constant. Our statement is nonasymptotic, and relies upon [29] and [6]. In particular, [29] requires $X$ and $\beta$ to satisfy the *Irrepresentable Condition*, which is sufficient to guarantee the exact recovery of the support of $\beta$ in some asymptotic regime; Section 3.3 connects with this line of work by showing that the *Irrepresentable Condition* holds with high probability under the stated assumptions.

As before, we have decided to state the theorem for a concrete value of $\lambda$, namely $2\sqrt{2\log p}$, but we could have used any value of the form $(1 + a)\sqrt{2\log p}$ with $a > 0$. When $a$ decreases, our proof indicates that one can lower the threshold on the minimum nonzero value of $\beta$ but that, at the same time, the probability of success is also lowered. When $a$ increases, the converse applies. Finally, our proof shows that, by setting $\lambda$ close to $\sqrt{2\log p}$ and imposing slightly stronger conditions on the coherence and the sparsity $S$, one can substantially lower the threshold on the minimum nonzero value of $\beta$ and bring it close to $\sigma\sqrt{2\log p}$.

We would also like to remark that, under the hypotheses of Theorem 1.3, one can somewhat improve the estimate (1.6) by using a two-step procedure similar to that proposed in [10]:

1. Use the lasso to find $\hat{I} \equiv \{i : \hat{\beta}_i \neq 0\}$.
2. Find $\tilde{\beta}$ by regressing $y$ onto the columns $(X_i)$, $i \in \hat{I}$.

Since $\hat{I} = I$ with high probability, we have that

$$\|X\tilde{\beta} - X\beta\|_{\ell_2}^2 = \|P[I]z\|_{\ell_2}^2$$

with high probability, where $P[I]$ is the projection onto the space spanned by the variables $(X_i)$. Because $\|P[I]z\|_{\ell_2}^2$ is concentrated around $|I| \cdot \sigma^2 = S \cdot \sigma^2$, it follows that, with high probability,

$$\|X\tilde{\beta} - X\beta\|_{\ell_2}^2 \leq C \cdot S \cdot \sigma^2,$$

where $C$ is a some small numerical constant. In other words, when the values of the nonzero entries of $\beta$ are sufficiently large, one does not have to pay the logarithmic factor.

1.4. *General model selection.* In many applications, $\beta$ is not sparse or does not have a real meaning, so that it does not make much sense to talk about the values of this vector. Consider an example to make this precise. Suppose we have noisy data $y$ (1.9) about an $n$-pixel digital image $f$, where $z$ is white noise. We wish to remove the noise (i.e., estimate the mean of the vector $y$). A majority of modern methods express the unknown signal as a superposition of fixed waveforms $(\varphi_i(t))_{1 \leq i \leq p}$,

$$(1.13) \qquad\qquad f(t) = \sum_{i=1}^{p} \beta_i \varphi_i(t),$$

and construct an estimate

$$\hat{f}(t) = \sum_{i=1}^{p} \hat{\beta}_i \varphi_i(t).$$

That is, one introduces a model $f = X\beta$, in which the columns of $X$ are the sampled waveforms $\varphi_i(t)$. It is now extremely popular to consider overcomplete representations with many more waveforms than samples (i.e., $p > n$). The reason for this is that overcomplete systems offer a wider range of generating elements that may be well suited to represent contributions from different phenomena; potentially, this wider range allows more flexibility in signal representation and enhances statistical estimation.

In this setup, two comments are in order. First, there is no ground truth associated with each coefficient $\beta_i$; there is no real wavelet or curvelet coefficient. Second, signals of general interest are never really exactly sparse; they are only approximately sparse, meaning that they may be well approximated by sparse expansions. These considerations emphasize the need to formulate results to cover those situations in which the precise values of $\beta_i$ are either ill-defined or meaningless.

In general, one can understand model selection as follows. Select a model— a subset $I$ of the columns of $X$—and construct an estimate of $X\beta$ by projecting $y$ onto the subspace generated by the variables in the model. Mathematically, this is formulated as

$$X\hat{\beta}[I] = P[I]y = P[I]X\beta + P[I]z,$$

where $P[I]$ denotes the projection onto the space spanned by the variables $(X_i)$, $i \in I$. What is the accuracy of $X\hat{\beta}[I]$? Note that

$$X\beta - X\hat{\beta}[I] = (\mathrm{Id} - P[I])X\beta - P[I]z;$$

therefore, the mean-squared error (MSE) obeys[3]

(1.14)        $$\mathbb{E}\|X\beta - X\hat{\beta}[I]\|^2 = \|(\mathrm{Id} - P[I])X\beta\|^2 + |I|\sigma^2.$$

This is the classical bias variance decomposition; the first term is the squared bias one gets by using only a subset of columns of $X$ to approximate the true vector $X\beta$. The second term is the variance of the estimator and is proportional to the size of the model $I$.

Hence, one can now define the *ideal model* achieving the minimum MSE over all models

(1.15)        $$\min_{I \subset \{1,\dots,p\}} \|(\mathrm{Id} - P[I])X\beta\|^2 + |I|\sigma^2.$$

See Figure 1 for a visual representation. We will refer to this as the ideal risk. It is ideal in the sense that one could achieve this performance if we had available an oracle which, knowing $X\beta$, would select for us the best model to use (i.e., the best subset of explanatory variables).

To connect this with our earlier discussion, one sees that, if there is a representation of $f = X\beta$ in which $\beta$ has $S$ nonzero terms, then the ideal risk is bounded by the variance term, namely, $S \cdot \sigma^2$ [just pick $I$ to be the support of $\beta$ in (1.15)]. The point we would like to make is that, whereas we did not search for an optimal bias-variance trade off in the previous section, we will here. The reason is that, even in the case where the model is interpretable, the projection estimate on the model corresponding to the nonzero values of $\beta_i$ may very well be inaccurate and have a mean-squared error that is far larger than (1.15). In particular, this typically happens if, out of the $S$ nonzero $\beta_i$'s, only a small fraction are really significant, while the others are not (e.g., in the sense that any individual test of significance would not reject the hypothesis that they vanish). In this sense, the main result of this section, Theorem 1.4, generalizes but also strengthens Theorem 1.2.

An important question is, of course, whether one can get close to the ideal risk (1.15) without the help of an oracle. It is known that solving the combinatorial optimization problem (1.2) with a value of $\lambda_0$ being a sufficiently large multiple of $\log p$ would provide an MSE within a multiplicative factor of order $\log p$ of the ideal risk. That real estimators with such properties exist is inspiring. Yet, solving (1.2) is computationally intractable. Our next result shows that, in a wide range problems, the lasso also nearly achieves the ideal risk.

---

[3]It is, again, simpler to state the performance in terms of expectation.

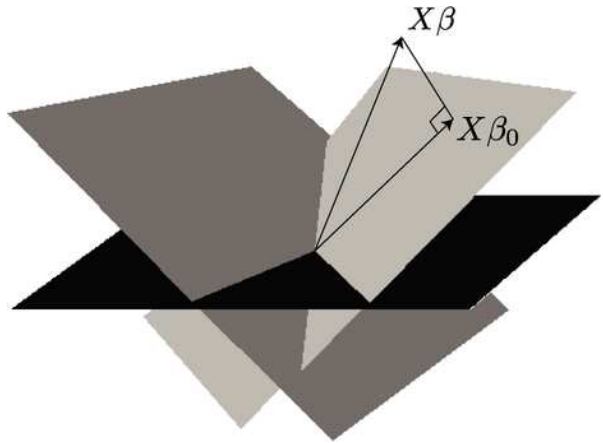

Fig. 1. *The vector $X\beta_0$ is the projection of $X\beta$ on an ideally selected subset of covariates. These covariates span a plane of optimal dimension, which, among all planes spanned by subsets of the same dimension, is closest to $X\beta$.*

We are naturally interested in quantifying the performance one can expect from the lasso in nearly all cases, and, just as before, we now introduce a useful statistical description of these cases. Consider the best model $I_0$ achieving the minimum in (1.15). In case of ties, pick one uniformly at random. Suppose $I_0$ is of cardinality $S$. Then, we introduce the *best $S$-dimensional subset model*, which is defined as follows:

1. The subset $I_0 \subset \{1, \dots, p\}$ of cardinality $S$ is distributed uniformly at random;

2. Define $\beta_0$ with support $I_0$ via

$$X\beta_0 = P[I_0]X\beta. \qquad (1.16)$$

   In other words, $\beta_0$ is the vector one would get by regressing the true mean vector $X\beta$ onto the variables in $I_0$; we call $\beta_0$ the ideal approximation. Conditional on $I_0$, the signs of the nonzero entries of $\beta_0$ are independent and equally likely to be $-1$ or $1$.

We make no assumption on the amplitudes. Our intent is just the same as before. All models are equally likely (there is no bias towards special variables), and one has no a priori information about the sign of the coefficients associated with each significant variable.

THEOREM 1.4. *Suppose that $X$ obeys the coherence property, and assume that the ideal approximation $\beta_0$ is taken from the best $S$-dimensional subset model. Suppose that $S \leq c_0 p/[\|X\|^2 \log p]$ for some positive numerical constant $c_0$. Then, the lasso estimate (1.3) computed with $\lambda = 2\sqrt{2 \log p}$*

*obeys*

$$(1.17) \quad \|X\beta - X\hat{\beta}\|_{\ell_2}^2$$
$$\leq (1 + \sqrt{2}) \left[ \min_{I \subset \{1, \ldots, p\}} \|X\beta - P[I]X\beta\|_{\ell_2}^2 + C_0' \, (2\log p) \cdot |I| \cdot \sigma^2 \right]$$

*with probability at least $1 - 6p^{-2\log 2} - p^{-1}(2\pi \log p)^{-1/2}$. The constant $C_0'$ may be taken as $12 + 10\sqrt{2}$.*

In words, the lasso nearly selects the best model in a very large majority of cases. This also strengthens our earlier result, since the right-hand side in (1.17) is always less or equal to $O(\log p)S\sigma^2$ whenever there is an $S$-sparse representation.[4]

Theorem 1.4 guarantees excellent performance in a broad range of problems. That is, if we have a design matrix $X$ whose columns are not too correlated, then, for most responses $X\beta$, the lasso will find a statistical model with low mean-squared error; simple extensions would also claim that the lasso finds a statistical model with very good predictive power, but we will not consider these here. As an illustrative example, we can consider predicting the clinical outcomes from different tumors on the basis of gene expression values for each of the tumors. In typical problems, one considers hundreds of tumors and tens of thousands of genes. While some of the gene expressions (the columns of $X$) are correlated, one can always eliminate redundant predictors (e.g., via clustering techniques). Once the statistician has designed an $X$ with low coherence, the lasso is guaranteed, in most cases, to find a subset of genes with near-optimal predictive power.

There is a slightly different formulation of this general result which may go as follows. Let $S_0$ be the maximum sparsity level $S_0 = \lfloor c_0 p / [\|X\|^2 \log p] \rfloor$, and, for each $S \leq S_0$, introduce $\mathcal{A}_S \subset \{-1, 0, 1\}^p$ as the set of all possible signs of vectors $\beta \in \mathbb{R}^p$, with $\text{sgn}(\beta_i) = 0$ if $\beta_i = 0$, such that exactly $S$ signs are nonzero. Then, under the hypotheses of our theorem, for each $X\beta \in \mathbb{R}^n$,

$$(1.18) \quad \|X\beta - X\hat{\beta}\|_{\ell_2}^2$$
$$\leq \min_{S \leq S_0} \min_{b \,:\, \text{sgn}(b) \in \mathcal{A}_{0,S}} (1 + \sqrt{2})[\|X\beta - Xb\|_{\ell_2}^2 + C_0' (2\log p) \cdot S \cdot \sigma^2]$$

---

[4] We have assumed that the mean response $f$ of interest is in the span of the columns of $X$ (i.e., of the form $X\beta$), which always happens when $p \geq n$ and $X$ has full column rank, for example. However, if this is not the case, the error would obey $\|f - X\hat{\beta}\|_{\ell_2}^2 = \|Pf - X\hat{\beta}\|_{\ell_2}^2 + \|(\text{Id} - P)f\|_{\ell_2}^2$, where $P$ is the projection onto the range of $X$. The first term obeys the oracle inequality, so that the lasso estimates $Pf$ in a near-optimal fashion. The second term is simply the size of the unmodelled part of the mean response.

with probability at least $1 - O(p^{-1})$, where one can still take $C'_0 = 12 + 10\sqrt{2}$ (the probability is with respect to the noise distribution). Above, $\mathcal{A}_{0,S}$ is a very large subset of $\mathcal{A}_S$, obeying

$$(1.19) \qquad |\mathcal{A}_{0,S}|/|\mathcal{A}_S| \geq 1 - O(p^{-1}).$$

Hence, for most $\beta$, the sub-oracle inequality (1.18) is actually the true oracle inequality.

For completeness, $\mathcal{A}_{0,S}$ is defined as follows. Let $b \in \mathcal{A}_S$ be supported on $I$; $b_I$ is the restriction of the vector $b$ to the index set $I$, and $X_I$ is the submatrix formed by selecting the columns of $X$ with indices in $I$. Then, we say that $b \in \mathcal{A}_{0,S}$ if and only if the following three conditions hold: (1) the submatrix $X_I^* X_I$ is invertible and obeys $\|(X_I^* X_I)^{-1}\| \leq 2$; (2) $\|X_{I^c}^* X_I (X_I^* X_I)^{-1} b_I\|_{\ell_\infty} \leq 1/4$ (recall that $b \in \{-1, 0, 1\}^p$ is a sign pattern); (3) $\max_{i \notin I} \|X_I (X_I^* X_I)^{-1} X_I^* X_i\| \leq c_0/\sqrt{\log p}$ for some numerical constant $c_0$. In Section 3, we will analyze these three conditions in detail and prove that $|\mathcal{A}_{0,S}|$ is large. The first condition is called the *invertibility condition*, and the second and third conditions are needed to establish that a certain *complementary size condition* holds (see Section 3).

1.5. *Implications for signal estimation.* Our findings may be of interest to researchers interested in signal estimation, and we now recast our main results in the language of signal processing. Suppose we are interested in estimating a signal $f(t)$ from observations

$$y(t) = f(t) + z(t), \qquad t = 0, \ldots, n - 1,$$

where $z$ is white noise with variance $\sigma^2$. We are given a dictionary of waveforms $(\varphi_i(t))_{1 \leq i \leq p}$, which are normalized so that $\sum_{t=0}^{n-1} \varphi_i^2(t) = 1$, and are looking for an estimate of the form $\hat{f}(t) = \sum_{i=1}^p \hat{\alpha}_i \varphi_i(t)$. When we have an overcomplete representation in which $p > n$, there are infinitely many ways of representing $f$ as a superposition of the dictionary elements.

Now, introduce the best $m$-term approximation $f_m$, which is defined via

$$\|f - f_m\|_{\ell_2} = \inf_{a : \#\{i, a_i \neq 0\} \leq m} \left\| f - \sum_i a_i \varphi_i \right\|_{\ell_2};$$

that is, it is that linear combination of at most $m$ elements of the dictionary that comes closest to the object $f$ of interest.[5] With these notations, if we could somehow guess the best model of dimension $m$, one would achieve an MSE equal to

$$\|f - f_m\|_{\ell_2}^2 + m\sigma^2.$$

---

[5] Note that, again, finding $f_m$ is generally a combinatorially hard problem.

Therefore, one can rewrite the ideal risk (which could be attained with the help of an oracle telling us exactly which subset of waveforms to use) as

$$\text{(1.20)} \qquad \min_{0 \le m \le p} \|f - f_m\|_{\ell_2}^2 + m\sigma^2,$$

which is exactly the trade-off between the approximation error and the number of terms in the partial expansion.[6]

Consider, now, the estimate $\hat{f} = \sum_i \hat{\alpha}_i \varphi_i$, where $\hat{\alpha}$ is the solution to

$$\text{(1.21)} \qquad \min_{a \in \mathbb{R}^p} \frac{1}{2} \left\| y - \sum_i a_i \varphi_i \right\|_{\ell_2}^2 + \lambda \sigma \|a\|_{\ell_1},$$

with $\lambda = 2\sqrt{2 \log p}$, say. Then, provided that the dictionary is not too redundant in the sense that $\max_{1 \le i < j \le p} |\langle \varphi_i, \varphi_j \rangle| \le c_0 / \log p$, Theorem 1.4 asserts that, for most signals $f$, the minimum-$\ell_1$ estimator (1.21) obeys

$$\text{(1.22)} \qquad \|\hat{f} - f\|_{\ell_2}^2 \le C_0 \left[ \inf_m \|f - f_m\|_{\ell_2}^2 + \log p \cdot m\sigma^2 \right]$$

with large probability and for some reasonably small numerical constant $C_0$. In other words, one obtains a squared error that is within a logarithmic factor of what can be achieved with information provided by a genie.

Overcomplete representations are now in widespread use, as in the field of artificial neural networks, for instance [12]. In computational harmonic analysis and image/signal processing, there is an emerging wisdom, which says that: (1) there is no universal representation for signals of interest, and (2) different representations are best for different phenomena ("best" is here understood as providing sparser representations). For instance:

- sinusoids are best for oscillatory phenomena;
- wavelets [18] are best for point-like singularities;
- curvelets [7, 8] are best for curve-like singularities (edges);
- local cosines are best for textures; and so on.

Thus, many efficient methods in modern signal estimation proceed by forming an overcomplete dictionary, a union of several distinct representations, and then extracting a sparse superposition that fits the data well. The main result of this paper says that, if one solves the quadratic program (1.21), then one is provably guaranteed near-optimal performance for most signals of interest, which is why these results might be of interest to people working in this field.

---

[6]It is also known that, for many interesting classes of signals $\mathcal{F}$ and appropriately chosen dictionaries, taking the supremum over $f \in \mathcal{F}$ in (1.20) comes within a log factor of the minimax risk for $\mathcal{F}$.

The spikes and sines model has been studied extensively in the literature on information theory in the nineties, and, there, the assumption that the "arrival times" of the spikes and the frequencies of the sinusoids are random is legitimate. In other situations, the model may be less adequate. For instance, in image processing, the large wavelet coefficients tend to appear early in the series, that is, at low frequencies. With this in mind, two comments are in order. First, it is likely that similar results would hold for other models (we just considered the simplest). And second, if we have a lot of a priori information about which coefficients are more likely to be significant, then we would probably not want to use the plain lasso (1.3) but rather incorporate this side information.

1.6. *Organization of the paper.* The paper is organized as follows. In Section 2, we explain why our results are nearly optimal and cannot be fundamentally improved. Section 3 introduces a recent result due to Joel Tropp, regarding the norm of certain random submatrices, which is essential to our proofs and proves all of our results. We conclude with a discussion in Section 4, where, for the most part, we relate our work with a series of other published results and distinguish our main contributions.

## 2. Optimality.

2.1. *For almost all sparse models.* A natural question is whether one can relax the condition about $\beta$ being *generically* sparse, or about $X\beta$ being well approximated by a *generically* sparse superposition of covariates. The emphasis is on "generic," meaning that our results apply to nearly all objects taken from a statistical ensemble but perhaps not all. This begs a question: can one hope to establish versions of our results that would hold universally? The answer is negative. Even in the case when $X$ has very low coherence, one can show that the lasso does not provide an accurate estimation of certain mean vectors $X\beta$ with a sparse coefficient sequence. This section gives one such example.

Suppose, as in Section 1.2, that we wish to estimate a signal assumed to be a sparse superposition of spikes and sinusoids. We assume that the length $n$ of the signal $f(t)$, $t = 0, 1, \ldots, n-1$, is equal to $n = 2^{2j}$ for some integer $j$. The basis of spikes is as before, and the orthobasis of sinusoids takes the form

$$\varphi_1(t) = 1/\sqrt{n},$$

$$\varphi_{2k}(t) = \sqrt{2/n}\cos(2\pi kt/n), \qquad k = 1, 2, \ldots, n/2 - 1,$$

$$\varphi_{2k+1}(t) = \sqrt{2/n}\sin(2\pi kt/n), \qquad k = 1, 2, \ldots, n/2 - 1,$$

$$\varphi_n(t) = (-1)^t/\sqrt{n}.$$

Recall the discrete identity (a discrete analog of the Poisson summation formula)

$$\sum_{k=0}^{2^j-1} \delta(t - k2^j) = \sum_{k=0}^{2^j-1} \frac{1}{\sqrt{n}} e^{i2\pi k2^j t/n}$$

$$\text{(2.1)} \qquad = \frac{1}{\sqrt{n}}(1 + (-1)^t) + \frac{2}{\sqrt{n}} \sum_{k=1}^{2^{j-1}-1} \cos(2\pi k2^j t/n)$$

$$= \varphi_1(t) + \varphi_n(t) + \sqrt{2} \sum_{k=1}^{2^{j-1}-1} \varphi_{k2^{j+1}}(t).$$

Then, consider the model

$$y = \mathbf{1} + z = X\beta + z,$$

where $\mathbf{1}$ is the constant signal equal to 1, and $X$ is the $n \times (2n-1)$ matrix

$$X = [\, I_n \quad F_{n,2:n} \,]$$

in which $I_n$ is the identity (the basis of spikes) and $F_{n,2:n}$ is the orthobasis of sinusoids minus the first basis vector $\varphi_1$. Note that this is a low-coherence matrix $X$, since $\mu(X) = \sqrt{2/n}$. In plain English, we are simply trying to estimate a constant-mean vector. It follows, from (2.1), that

$$1 = \sqrt{n} \left[ \sum_{k=0}^{2^j-1} \delta(t - k2^j) - \varphi_n(t) - \sqrt{2} \sum_{k=1}^{2^{j-1}-1} \varphi_{k2^{j+1}}(t) \right],$$

so that $\mathbf{1}$ has a sparse expansion, since it is a superposition of at most $\sqrt{n}$ spikes and $\sqrt{n}/2$ sinusoids (it can also be deduced from existing results that this is actually the sparsest expansion). In other words, if we knew which column vectors to use, one could obtain

$$\mathbb{E}\|X\beta^\star - X\beta\|_{\ell_2}^2 = \tfrac{3}{2}\sqrt{n}\sigma^2.$$

How does the lasso compare? We claim that, with very high probability,

$$\text{(2.2)} \qquad \hat{\beta}_i = \begin{cases} y_i - \lambda\sigma, & i \in \{1, \ldots, n\}, \\ 0, & i \in \{n+1, \ldots, 2n-1\}, \end{cases}$$

so that

$$\text{(2.3)} \qquad X\hat{\beta} = y - \lambda\sigma\mathbf{1},$$

provided that $\lambda\sigma \leq 1/2$. In short, the lasso does not find the sparsest model at all. As a matter of fact, it finds a model as dense as it can be, and the resulting mean-squared error is awful, since

$$\mathbb{E}\|X\hat{\beta} - X\beta\|_{\ell_2}^2 \approx (1 + \lambda^2)n\sigma^2.$$

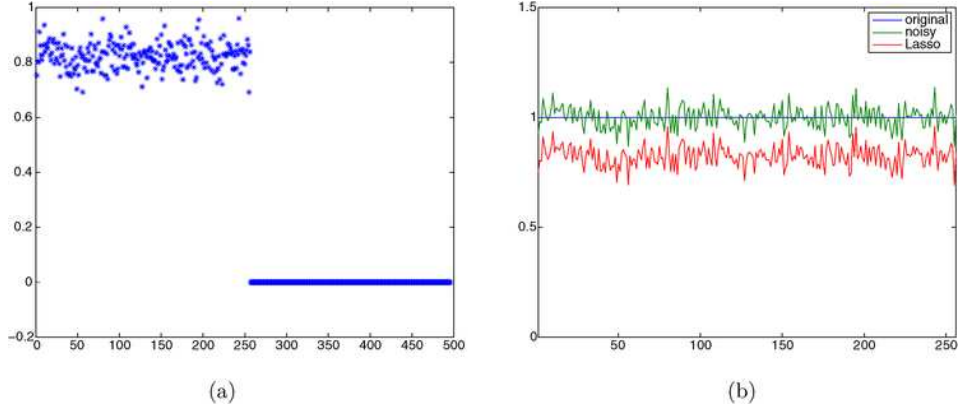

Fig. 2.  *Sparse signal recovery with the lasso.* (a) *Values of the estimated coefficients. All the spike coefficients are obtained by soft-thresholding* $y$ *and are nonzero.* (b) *Lasso signal estimate;* $X\hat{\beta}$ *is just a shifted version of the noisy signal.*

Even if one could somehow remove the bias, this would still be a very bad performance.

An illustrative numerical example is displayed in Figure 2. In this example, $n = 256$ so that $p = 512 - 1 = 511$. The mean vector $X\beta$ is made up as above, and there is a representation in which $\beta$ has only 24 nonzero coefficients. Yet, the lasso finds a model of dimension 256 (i.e., select as many variables as there are observations).

We need to justify (2.2), as (2.3) would be an immediate consequence. It follows, from taking the subgradient of the lasso functional, that $\hat{\beta}$ is a minimizer if and only if

$$(2.4) \qquad \begin{aligned} X_i^*(y - X\hat{\beta}) &= \lambda\sigma\operatorname{sgn}(\hat{\beta}_i), \qquad \hat{\beta}_i \neq 0, \\ |X_i^*(y - X\hat{\beta})| &\leq \lambda\sigma, \qquad \hat{\beta}_i = 0. \end{aligned}$$

One can further establish that $\hat{\beta}$ is the unique minimizer of (1.3) if

$$(2.5) \qquad \begin{aligned} X_i^*(y - X\hat{\beta}) &= \lambda\sigma\operatorname{sgn}(\hat{\beta}_i), \qquad \hat{\beta}_i \neq 0, \\ |X_i^*(y - X\hat{\beta})| &< \lambda\sigma, \qquad \hat{\beta}_i = 0, \end{aligned}$$

and the columns indexed by the support of $\hat{\beta}$ are linearly independent (note the strict inequalities). We then simply need to show that $\hat{\beta}$, given by (2.2), obeys (2.5). Suppose that $\min_i y_i > \lambda\sigma$. A sufficient condition is that $\max_i |z_i| < 1 - \lambda\sigma$, which occurs with very large probability if $\lambda\sigma \leq 1/2$ and $\lambda > \sqrt{2\log n}$ [see (3.4) with $X = I$]. (One can always allow for larger noise by multiplying the signal by a factor greater than 1.) Note that $y - X\hat{\beta} = \lambda\sigma\mathbf{1}$, so that, for $i \in \{1, \ldots, n\}$, we have

$$X_i^*(y - X\hat{\beta}) = \lambda \sigma = \lambda \sigma \, \text{sgn}(\hat{\beta}_i),$$

whereas, for $i \in \{n+1, \ldots, 2n-1\}$, we have

$$X_i^*(y - X\hat{\beta}) = \lambda \sigma \langle X_i, \mathbf{1} \rangle = 0,$$

which proves our claim.

To summarize, even when the coherence is low (i.e., of size about $1/\sqrt{n}$) there are sparse vectors $\beta$ with sparsity level about equal to $\sqrt{n}$, for which the lasso completely misbehaves (we presented an example but there are of course many others). Therefore, it is a fact that none of our theorems, namely, Theorems 1.2–1.4, can hold for all $\beta$'s. In this sense, they are sharp.

2.2. *For sufficiently incoherent matrices.* We now show that predictors cannot be too collinear, and we begin by examining a small problem in which $X$ is a $2 \times 2$ matrix $X = [X_1, X_2]$. We violate the coherence property by choosing $X_1$ and $X_2$, so that $\langle X_1, X_2 \rangle = 1 - \varepsilon$, where we think of $\varepsilon$ as being very small. Assume, without loss of generality, that $\sigma = 1$ to simplify. Now, consider

$$\beta = \frac{a}{\varepsilon} \begin{bmatrix} 1 \\ -1 \end{bmatrix},$$

where $a$ is some positive amplitude and observe that $X\beta = a\varepsilon^{-1}(X_1 - X_2)$, and $X^*X\beta = a(1, -1)^*$. It is well known that the lasso estimate $\hat{\beta}$ vanishes if $\|X^*y\|_{\ell_\infty} \leq \lambda$. Now,

$$\|X^*y\|_{\ell_\infty} \leq a + \|X^*z\|_{\ell_\infty},$$

so that, if $a = 1$, say, and $\lambda$ is not ridiculously small, then there is a positive probability $\pi_0$ that $\hat{\beta} = 0$, where $\pi_0 \geq \mathbb{P}(\|X^*z\|_\infty \leq \lambda - 1)$.[7] For example, if $\lambda > 1 + 3 = 4$, then $\hat{\beta} = 0$, as long as both entries of $X^*z$ are within 3 standard deviations of 0. When $\hat{\beta} = 0$, the squared error loss obeys

$$\|X\beta\|_{\ell_2}^2 = 2\frac{a^2}{\varepsilon},$$

which can be made arbitrarily large if we allow $\varepsilon$ to be arbitrarily small.

Of course, the culprit in our 2-by-2 example is hardly sparse, and we now consider the $n \times n$ diagonal block matrix $X_0$ ($n$ even)

$$X_0 = \begin{bmatrix} X & & & \\ & X & & \\ & & \ddots & \\ & & & X \end{bmatrix}$$

---

[7] $\pi_0$ can be calculated since $X^*z$ is a bivariate Gaussian variable.

with blocks made out of $n/2$ copies of $X$. We now consider $\beta$ from the $S$-sparse model with independent entries sampled from the distribution (we choose $a = 1$ for simplicity but we could consider other values as well)

$$\beta_i = \begin{cases} \varepsilon^{-1}, & \text{w. p.} \ n^{-1/2}, \\ -\varepsilon^{-1}, & \text{w. p.} \ n^{-1/2}, \\ 0, & \text{w. p.} \ 1 - 2n^{-1/2}. \end{cases}$$

Certainly, the support of $\beta$ is random, and the signs are random. One could argue that the size of the support is not fixed (the expected value is $2\sqrt{n}$, so that $\beta$ is sparse with very large probability) but this is obviously unessential.[8]

Because $X_0$ is block diagonal, the lasso functional becomes additive, and the lasso will minimize each individual term of the form $\frac{1}{2}\|Xb^{(i)} - y^{(i)}\|_{\ell_2}^2 + \lambda\|b^{(i)}\|_{\ell_1}$, where $b^{(i)} = (b_{2i-1}, b_{2i})$ and $y^{(i)} = (y_{2i-1}, y_{2i})$. If, for any of these subproblems, $\beta^{(i)} = \pm\varepsilon^{-1}(1, -1)$ as in our 2-by-2 example above, then the squared error will blow up (as $\varepsilon$ gets smaller) with probability $\pi_0$. With $i$ fixed, $\mathbb{P}(\beta^{(i)} = \pm\varepsilon^{-1}(1, -1)) = 2/n$ and, thus, the probability that none of these sub-problems is poised to blow up is $(1 - \frac{2}{n})^{n/2} \to \frac{1}{e}$ as $n \to \infty$. Formalizing matters, we have a squared loss of at least $2/\varepsilon$ with probability at least $\pi_0(1 - (1 - \frac{2}{n})^{n/2})$. Note that, when $n$ is large, $\lambda$ is large, so that $\pi_0$ is close to 1, and the lasso badly misbehaves with a probability greater or equal to a quantity approaching $1 - 1/e$.

In conclusion, the lasso may perform badly, even with a random $\beta$, when all our assumptions are met but the coherence property. To summarize, an upper bound on the coherence is also necessary.

**3. Proofs.** In this section, we prove all of our results. It is sufficient to establish our theorems with $\sigma = 1$, as the general case is treated by a simple rescaling. Therefore, we conveniently assume that $\sigma = 1$ from now on. Here, and in the remainder of this paper, $x_I$ is the restriction of the vector $x$ to an index set $I$, and, for a matrix $X$, $X_I$ is the submatrix formed by selecting the columns of $X$ with indices in $I$. In the following, it will also be convenient to denote, by $K$, the functional

$$(3.1) \qquad\qquad K(y, b) = \tfrac{1}{2}\|y - Xb\|_{\ell_2}^2 + 2\lambda_p\|b\|_{\ell_1},$$

in which $\lambda_p = \sqrt{2 \log p}$.

---

[8]We could alternatively select the support at random and randomly assign the signs, and this would not change our story in the least.

3.1. *Preliminaries.* We will make frequent use of subgradients, and we begin by briefly recalling what these are. We say that $u \in \mathbb{R}^p$ is a subgradient of a convex function $f \colon \mathbb{R}^p \to \mathbb{R}$ at $x_0$ if $f$ obeys

$$(3.2) \qquad f(x) \geq f(x_0) + \langle u, x - x_0 \rangle$$

for all $x$.

Further, our arguments will repeatedly use two general results that we now record. The first states that the lasso estimate is feasible for the Dantzig selector optimization problem.

LEMMA 3.1. *The lasso estimate obeys*

$$(3.3) \qquad \|X^*(y - X\hat{\beta})\|_{\ell_\infty} \leq 2\lambda_p.$$

PROOF. Since $\hat{\beta}$ minimizes $f(b) = K(y, b)$ over $b$, 0 must be a subgradient of $f$ at $\hat{\beta}$. Now, the subgradients of $f$ at $b$ are of the form

$$X^*(Xb - y) + 2\lambda_p \varepsilon,$$

where $\varepsilon$ is any $p$-dimensional vector obeying $\varepsilon_i = \mathrm{sgn}(b_i)$ if $b_i \neq 0$ and $|\varepsilon_i| \leq 1$ otherwise. Hence, since 0 is a subgradient at $\hat{\beta}$, there exists $\varepsilon$ as above such that

$$X^*(X\hat{\beta} - y) = -2\lambda_p \varepsilon.$$

The conclusion follows from $\|\varepsilon\|_{\ell_\infty} \leq 1$. $\square$

The second general result states that $\|X^*z\|_{\ell_\infty}$ cannot be too large. With large probability, $z \sim \mathcal{N}(0, I)$ obeys

$$(3.4) \qquad \|X^*z\|_{\ell_\infty} = \max_i |\langle X_i, z \rangle| \leq \lambda_p.$$

This is standard and simply follows from the fact that $\langle X_i, z \rangle \sim \mathcal{N}(0, 1)$. Hence, for each $t > 0$,

$$(3.5) \qquad \mathbb{P}(\|X^*z\|_{\ell_\infty} > t) \leq 2p \cdot \phi(t)/t,$$

where $\phi(t) \equiv (2\pi)^{-1/2} e^{-t^2/2}$. Better bounds may be possible, but we will not pursue these refinements here. Also, note that $\|X^*z\|_{\ell_\infty} \leq \sqrt{2}\lambda_p$ with probability at least $1 - p^{-1}(2\pi \log p)^{-1/2}$. These two general facts have an interesting consequence, since it follows from the decomposition $y = X\beta + z$ and the triangle inequality that, with high probability,

$$(3.6) \qquad \begin{aligned} \|X^*X(\beta - \hat{\beta})\|_{\ell_\infty} &\leq \|X^*(X\beta - y)\|_{\ell_\infty} + \|X^*(y - X\hat{\beta})\|_{\ell_\infty} \\ &= \|X^*z\|_{\ell_\infty} + \|X^*(y - X\hat{\beta})\|_{\ell_\infty} \\ &\leq (\sqrt{2} + 2)\lambda_p. \end{aligned}$$

3.2. *Proof of Theorem 1.2.* Put $I$ for the support of $\beta$. To prove our claim, we first establish that (1.6) holds provided that the following three deterministic conditions are satisfied:

- *Invertibility condition.* The submatrix $X_I^* X_I$ is invertible and obeys

$$(3.7) \qquad \|(X_I^* X_I)^{-1}\| \leq 2.$$

  The number 2 is arbitrary; we just need the smallest eigenvalue of $X_I^* X_I$ to be bounded away from zero.
- *Orthogonality condition.* The vector $z$ obeys $\|X^* z\|_{\ell_\infty} \leq \sqrt{2}\lambda_p$.
- *Complementary size condition.* The following inequality holds

$$(3.8) \qquad \begin{aligned} \|X_{I^c}^* X_I (X_I^* X_I)^{-1} X_I^* z\|_{\ell_\infty} &+ 2\lambda_p \|X_{I^c}^* X_I (X_I^* X_I)^{-1} \mathrm{sgn}(\beta_I)\|_{\ell_\infty} \\ &\leq (2 - \sqrt{2})\lambda_p. \end{aligned}$$

This section establishes the main estimate (1.6), assuming these three conditions hold, whereas the next will show that all three conditions hold with large probability, hence proving Theorem 1.2. Note that, when $z$ is white noise, we already know that the orthogonality condition holds with probability at least $1 - p^{-1}(2\pi \log p)^{-1/2}$.

Assume, then, that all three conditions above hold. Since $\hat{\beta}$ minimizes $K(y, b)$, we have $K(y, \hat{\beta}) \leq K(y, \beta)$ or, equivalently,

$$\tfrac{1}{2}\|y - X\hat{\beta}\|_{\ell_2}^2 + 2\lambda_p \|\hat{\beta}\|_{\ell_1} \leq \tfrac{1}{2}\|y - X\beta\|_{\ell_2}^2 + 2\lambda_p \|\beta\|_{\ell_1}.$$

Set $h = \hat{\beta} - \beta$, and note that

$$\|y - X\hat{\beta}\|_{\ell_2}^2 = \|(y - X\beta) - Xh\|_{\ell_2}^2 = \|Xh\|_{\ell_2}^2 + \|y - X\beta\|_{\ell_2}^2 - 2\langle Xh, y - X\beta\rangle.$$

Plugging this identity with $z = y - X\beta$ into the above inequality and rearranging the terms gives

$$(3.9) \qquad \tfrac{1}{2}\|Xh\|_{\ell_2}^2 \leq \langle Xh, z\rangle + 2\lambda_p(\|\beta\|_{\ell_1} - \|\hat{\beta}\|_{\ell_1}).$$

Next, break $h$ up into $h_I$ and $h_{I^c}$ (observe that $\hat{\beta}_{I^c} = h_{I^c}$) and rewrite (3.9) as

$$\tfrac{1}{2}\|Xh\|_{\ell_2}^2 \leq \langle h, X^* z\rangle + 2\lambda_p(\|\beta_I\|_{\ell_1} - \|\beta_I + h_I\|_{\ell_1} - \|h_{I^c}\|_{\ell_1}).$$

For each $i \in I$, we have

$$|\hat{\beta}_i| = |\beta_i + h_i| \geq |\beta_i| + \mathrm{sgn}(\beta_i)\, h_i$$

and, thus, $\|\beta_I + h_I\|_{\ell_1} \geq \|\beta\|_{\ell_1} + \langle h_I, \mathrm{sgn}(\beta_I)\rangle$. Inserting this inequality above yields

$$(3.10) \qquad \tfrac{1}{2}\|Xh\|_{\ell_2}^2 \leq \langle h, X^* z\rangle - 2\lambda_p(\langle h_I, \mathrm{sgn}(\beta_I)\rangle + \|h_{I^c}\|_{\ell_1}).$$

Observe, now, that $\langle h, X^*z \rangle = \langle h_I, X_I^*z \rangle + \langle h_{I^c}, X_{I^c}^*z \rangle$ and that the orthogonality condition implies

$$\langle h_{I^c}, X_{I^c}^*z \rangle \le \|h_{I^c}\|_{\ell_1}\|X_{I^c}^*z\|_{\ell_\infty} \le \sqrt{2}\lambda_p \|h_{I^c}\|_{\ell_1}.$$

The conclusion is the useful estimate

$$(3.11) \qquad \tfrac{1}{2}\|Xh\|_{\ell_2}^2 \le \langle h_I, v \rangle - (2-\sqrt{2})\lambda_p\|h_{I^c}\|_{\ell_1},$$

where $v \equiv X_I^*z - 2\lambda_p \operatorname{sgn}(\beta_I)$.

We complete the argument by bounding $\langle h_I, v \rangle$. The key here is to use the fact that $\|X^*Xh\|_{\ell_\infty}$ is known to be small, as pointed out by Terence Tao [25]. We have

$$(3.12) \qquad \begin{aligned}
\langle h_I, v \rangle &= \langle (X_I^*X_I)^{-1}X_I^*X_Ih_I, v \rangle \\
&= \langle X_I^*X_Ih_I, (X_I^*X_I)^{-1}v \rangle \\
&= \langle X_I^*Xh, (X_I^*X_I)^{-1}v \rangle - \langle X_I^*X_{I^c}h_{I^c}, (X_I^*X_I)^{-1}v \rangle \equiv A_1 - A_2.
\end{aligned}$$

We address each of the two terms individually. First,

$$A_1 \le \|X_I^*Xh\|_{\ell_\infty} \cdot \|(X_I^*X_I)^{-1}v\|_{\ell_1}$$

and

$$\begin{aligned}
\|(X_I^*X_I)^{-1}v\|_{\ell_1} &\le \sqrt{S} \cdot \|(X_I^*X_I)^{-1}v\|_{\ell_2} \\
&\le \sqrt{S} \cdot \|(X_I^*X_I)^{-1}\|\|v\|_{\ell_2} \\
&\le S \cdot \|(X_I^*X_I)^{-1}\|\|v\|_{\ell_\infty}.
\end{aligned}$$

Consider the following: (1) $\|X_I^*Xh\|_{\ell_\infty} \le (2+\sqrt{2})\,\lambda_p$ by Lemma 3.1 together with the orthogonality condition [see (3.6)], and (2) $\|(X_I^*X_I)^{-1}\|_{\ell_2} \le 2$ by the invertibility condition. Because of this, we have

$$A_1 \le 2(2+\sqrt{2})\lambda_p S\|v\|_{\ell_\infty}.$$

However,

$$\|v\|_{\ell_\infty} \le \|X_I^*z\|_{\ell_\infty} + 2\lambda_p \le (2+\sqrt{2})\lambda_p,$$

so that

$$(3.13) \qquad A_1 \le 2(2+\sqrt{2})^2\lambda_p^2 \cdot S.$$

Second, we simply bound the other term $A_2 = \langle h_{I^c}, X_{I^c}^*X_I(X_I^*X_I)^{-1}v \rangle$ by

$$|A_2| \le \|h_{I^c}\|_{\ell_1}\|X_{I^c}^*X_I(X_I^*X_I)^{-1}v\|_{\ell_\infty}$$

with $v = X_I^*z - 2\lambda_p \operatorname{sgn}(\beta_I)$. Since

$$\begin{aligned}
\|X_{I^c}^*X_I&(X_I^*X_I)^{-1}v\|_{\ell_\infty} \\
&\le \|X_{I^c}^*X_I(X_I^*X_I)^{-1}X_I^*z\|_{\ell_\infty} + 2\lambda_p\|X_{I^c}^*X_I(X_I^*X_I)^{-1}\operatorname{sgn}(\beta_T)\|_{\ell_\infty} \\
&\le (2-\sqrt{2})\lambda_p
\end{aligned}$$

because of the complementary size condition, we have

$$|A_2| \leq (2 - \sqrt{2}) \lambda_p \|h_{I^c}\|_{\ell_1}.$$

To summarize,

$$(3.14) \qquad |\langle h_I, v \rangle| \leq 2(2 + \sqrt{2})^2 \lambda_p^2 \cdot S + (2 - \sqrt{2}) \lambda_p \|h_{I^c}\|_{\ell_1}.$$

We conclude by inserting (3.14) into (3.11), which gives

$$\tfrac{1}{2} \|X(\hat\beta - \beta)\|_{\ell_2}^2 \leq 2(2 + \sqrt{2})^2 \lambda_p^2 \cdot S,$$

which is what we needed to prove.

3.3. *Norms of random submatrices.* In this section, we establish that the invertibility and the complementary size conditions hold with large probability. These essentially rely on a recent result of Joel Tropp, which we state first.

THEOREM 3.2 [27]. *Suppose that a set $I$ of predictors is sampled using a Bernoulli model by first creating a sequence $(\delta_j)_{1 \leq j \leq p}$ of i.i.d. random variables with $\delta_j = 1$ w.p. $S/p$ and $\delta_j = 0$ w.p. $1 - S/p$, and then setting $I \equiv \{j : \delta_j = 1\}$ so that $\mathbb{E}|I| = S$. Then, for $q = 2 \log p$,*

$$(3.15) \qquad (\mathbb{E}\|X_I^* X_I - \mathrm{Id}\,\|^q)^{1/q} \leq 30 \mu(X) \log p + 13 \sqrt{\frac{2S\|X\|^2 \log p}{p}}$$

*provided that $S\|X\|^2/p \leq 1/4$. In addition, for the same value of $q$*

$$(3.16) \qquad \left(\mathbb{E} \max_{i \in I^c} \|X_I^* X_i\|_{\ell^2}^q\right)^{1/q} \leq 4 \mu(X) \sqrt{\log p} + \sqrt{S\|X\|^2/p}.$$

The first inequality (3.15) can be derived from the last equation in Section 4 of [27]. To be sure, using the notations of that paper and letting $H \equiv X^* X - \mathrm{Id}$, Tropp shows that

$$\mathbb{E}_q\|RHR\| \leq 15\bar q \mathbb{E}_q\|RHR'\|_{\max} + 12\sqrt{\delta\bar q}\|HR\|_{1 \to 2} + 2\delta\|H\|, \qquad \delta = S/p,$$

where $\bar q = \max\{q, 2 \log p\}$. Now, consider the following three facts: (1) $\|R \times HR'\|_{\max} \leq \mu(X)$; (2) $\|HR\|_{1 \to 2} \leq \|X\|$; and (3) $\|H\| \leq \|X\|^2$. The first assertion is immediate. The second is justified in [27]. For the third, observe that $\|X^* X - \mathrm{Id}\| \leq \max\{\|X\|^2 - 1, 1\}$ (this is an equality when $p > n$), and the claim follows from $\|X\| \geq 1$, which holds, since $X$ has unit-normed columns. With $q = 2 \log p$, this gives

$$\mathbb{E}_q\|RHR\| \leq 30 \mu(X) \log p + 12 \sqrt{\frac{2S \log p \|X\|^2}{p}} + \frac{2S\|X\|^2}{p}.$$

Suppose that $S\|X\|^2/p \leq 1/4$; then, we can simplify the above inequality and obtain

$$\mathbb{E}_q\|RHR\| \leq 30\mu(X)\log p + (12\sqrt{2\log p}+1)\sqrt{S\|X\|^2/p},$$

which implies (3.15). The second inequality (3.16) is exactly Corollary 5.1 in [27].

The inequalities (3.15) and (3.16) also hold for our slightly different model, in which $I \subset \{1,\ldots,p\}$ is a random subset of predictors with $S$ elements, provided that the right-hand side of both inequalities be multiplied by $2^{1/q}$. This follows from a simple Poissonization argument, which is similar to that posed in the proof of Lemma 3.6.

It is now time to investigate how these results imply our conditions, and we first examine how (3.15) implies the invertibility condition. Let $I$ be a random set and put $Z = \|X_I^*X_I - \mathrm{Id}\|$. Clearly, if $Z \leq 1/2$, then all the eigenvalues of $X_I^*X_I$ are in the interval $[1/2, 3/2]$ and $\|(X_I^*X_I)^{-1}\| \leq 2$. Suppose that $\mu(X)$ and $S$ are sufficiently small, so that the right-hand side of (3.15) is less than $1/4$, say. This happens when the coherence $\mu(X)$ and $S$ obey the hypotheses of the theorem. Then, by Markov's inequality, we have that, for $q = 2\log p$,

$$\mathbb{P}(Z > 1/2) \leq 2^q\mathbb{E}Z^q \leq (1/2)^q.$$

In other words, the invertibility condition holds with probability exceeding $1 - p^{-2\log 2}$.

Recalling that the signs of the nonzero entries of $\beta$ are i.i.d. symmetric variables, we now examine the complementary size condition and begin with a simple lemma.

LEMMA 3.3. Let $(W_j)_{j\in J}$ be a fixed collection of vectors in $\ell_2(I)$ and consider the random variable $Z_0$ defined by $Z_0 = \max_{j\in J} |\langle W_j, \mathrm{sgn}(\beta_I)\rangle|$. Then,

$$(3.17) \qquad \mathbb{P}(Z_0 \geq t) \leq 2|J| \cdot e^{-t^2/2\kappa^2}$$

for any $\kappa$ obeying $\kappa \geq \max_{j\in J} \|W_j\|_{\ell_2}$. Similarly, letting $(W_j')_{j\in J}$ be a fixed collection of vectors in $\mathbb{R}^n$ and setting $Z_1 = \max_{j\in J} |\langle W_j', z\rangle|$, we have

$$(3.18) \qquad \mathbb{P}(Z_1 \geq t) \leq 2|J| \cdot e^{-t^2/2\kappa^2}$$

for any $\kappa$ obeying $\kappa \geq \max_{j\in J} \|W_j'\|_{\ell_2}$.[9]

---

[9] Note that this lemma also holds if the collection of vectors $(W_j)_{j\in J}$ is random, as long as it is independent of $\mathrm{sgn}(\beta_I)$ and $z$.

PROOF.  The first inequality is an application of Hoeffding's inequality. Indeed, letting $Z_{0,j} = \langle W_j, \operatorname{sgn}(\beta_I) \rangle$, Hoeffding's inequality gives

$$(3.19) \qquad \mathbb{P}(|Z_{0,j}| > t) \le 2e^{-t^2/2\|W_j\|_{\ell_2}^2} \le 2e^{-t^2/2\max_j \|W_j\|_{\ell_2}^2}.$$

Inequality (3.17) then follows from the union bound. The second part is even easier, since $Z_{1,j} = \langle W'_j, z \rangle \sim \mathcal{N}(0, \|W'_j\|_{\ell_2}^2)$; thus,

$$(3.20) \qquad \mathbb{P}(|Z_{1,j}| > t) \le 2e^{-t^2/2\|W'_j\|_{\ell_2}^2} \le 2e^{-t^2/2\max_j \|W'_j\|_{\ell_2}^2}.$$

Again, the union bound gives (3.18).  $\square$

For each $i \in I^c$, define $Z_{0,i}$ and $Z_{1,i}$ as

$$Z_{0,i} = X_i^* X_I (X_I^* X_I)^{-1} \operatorname{sgn}(\beta_I) \quad \text{and} \quad Z_{1,i} = X_i^* X_I (X_I^* X_I)^{-1} X_I^* z.$$

With these notations, in order to prove the complementary size condition, it is sufficient to show that, with large probability,

$$2\lambda_p Z_0 + Z_1 \le (2 - \sqrt{2})\lambda_p,$$

where $Z_0 = \max_{i \in I^c} |Z_{0,i}|$ and likewise for $Z_1$. Therefore, it is sufficient to prove that, with large probability,

$$Z_0 \le 1/4 \quad \text{and} \quad Z_1 \le (3/2 - \sqrt{2})\lambda_p.$$

The idea is of course to apply Lemma 3.3 together with Theorem 3.2. We have

$$Z_{0,i} = \langle W_i, \operatorname{sgn}(\beta_I) \rangle \quad \text{and} \quad Z_{1,i} = \langle W'_i, z \rangle,$$

where

$$W_i = (X_I^* X_I)^{-1} X_I^* X_i \quad \text{and} \quad W'_i = X_I (X_I^* X_I)^{-1} X_I^* X_i.$$

Recall the definition of $Z$ above and consider the event $E = \{Z \le 1/2\} \cap \{\max_{i \in I^c} \|X_I^* X_i\| \le \gamma\}$ for some positive $\gamma$. On this event, all the singular values of $X_I$ are between $1/\sqrt{2}$ and $\sqrt{3/2}$; thus, $\|(X_I^* X_I)^{-1}\| \le 2$ and $\|X_I (X_I^* X_I)^{-1}\| \le \sqrt{2}$, which gives

$$\|W_i\| \le 2\gamma \quad \text{and} \quad \|W'_i\| \le \sqrt{2}\gamma.$$

Applying (3.17) and (3.18) gives

$$\begin{aligned}
\mathbb{P}(\{Z_0 \ge t\} \cup \{Z_1 \ge u\}) &\le \mathbb{P}(\{Z_0 \ge t\} \cup \{Z_1 \ge u\} | E) + \mathbb{P}(E^c) \\
&\le \mathbb{P}(Z_0 \ge t | E) + \mathbb{P}(Z_1 \ge u \mid E) + \mathbb{P}(E^c) \\
&\le 2p\, e^{-t^2/8\gamma^2} + 2p e^{-u^2/4\gamma^2} \\
&\quad + \mathbb{P}(Z > 1/2) + \mathbb{P}\left(\max_{i \in I^c} \|X_I^* X_i\| > \gamma\right).
\end{aligned}$$

We already know that the second to last term of the right-hand side is less than $p^{-2\log 2}$, provided that $\mu(X)$ and $S$ obey the conditions of the theorem. For the other three terms, let $\gamma_0$ be the right-hand side of (3.16). For $t = 1/4$, one can find a constant $c_0$ such that, if $\gamma < c_0/\sqrt{\log p}$, then $2pe^{-t^2/8\gamma^2} \leq 2p^{-2\log 2}$, say. Likewise, for $u = (3/2 - \sqrt{2})\lambda_p$, we may have $2pe^{-u^2/4\gamma^2} \leq 2p^{-2\log 2}$. The last term is treated by Markov's inequality, since, for $q = 2\log p$, (3.16) gives

$$\mathbb{P}\Big(\max_{i\in I^c} \|X_I^* X_i\| > \gamma\Big) \leq \gamma^{-q} \cdot \mathbb{E}\Big(\max_{i\in I^c} \|X_I^* X_i\|^q\Big) \leq (\gamma_0/\gamma)^q.$$

Therefore, if $\gamma_0 \leq \gamma/2 = c_0/2\sqrt{\log p}$, we have that this last term does not exceed $1 - p^{-2\log 2}$. For $\mu(X)$ and $S$ obeying the hypotheses of Theorem 1.2, it is indeed the case that $\gamma_0 \leq c_0/2\sqrt{\log p}$. In conclusion, we have shown that all three conditions hold under our hypotheses with probability at least $1 - 6p^{-2\log 2} - p^{-1}(2\pi\log p)^{-1/2}$.

In passing, we would like to remark that proving that $Z_0 \leq 1/4$ establishes that the strong irrepresentable condition from [29] holds (with high probability). This condition states that, if $I$ is the support of $\beta$,

$$\|X_{I^c}^* X_I (X_I^* X_I)^{-1}\mathrm{sgn}(\beta_I)\|_{\ell_\infty} \leq 1 - \nu,$$

where $\nu$ is any (small) constant greater than zero (this condition is used to show the asymptotic recovery of the support of $\beta$).

3.4. *Proof of Theorem 1.4.* The proof of Theorem 1.4 parallels that of Theorem 1.2, and we only sketch it, although we carefully detail the main differences. Let $I_0$ be the support of $\beta_0$. Just as before, all three conditions of Section 3.2, with $I_0$ in place of $I$ and $\beta_0$ in place of $\beta$, hold with overwhelming probability. From now on, we just assume that they are all true.

Since $\hat\beta$ minimizes $K(y, b)$, we have $K(y, \hat\beta) \leq K(y, \beta_0)$ or, equivalently,

$$(3.21) \qquad \tfrac{1}{2}\|y - X\hat\beta\|_{\ell_2}^2 + 2\lambda_p\|\hat\beta\|_{\ell_1} \leq \tfrac{1}{2}\|y - X\beta_0\|_{\ell_2}^2 + 2\lambda_p\|\beta_0\|_{\ell_1}.$$

Expand $\|y - X\hat\beta\|_{\ell_2}^2$ as

$$\|y - X\hat\beta\|_{\ell_2}^2 = \|z - (X\hat\beta - X\beta)\|_{\ell_2}^2 = \|z\|_{\ell_2}^2 - 2\langle z, X\hat\beta - X\beta\rangle + \|X\hat\beta - X\beta\|_{\ell_2}^2$$

and $\|y - X\beta_0\|_{\ell_2}^2$ in the same way. Then, plug these identities into (3.21) to obtain

$$(3.22) \qquad \begin{aligned} \tfrac{1}{2}\|X\hat\beta - X\beta\|_{\ell_2}^2 &\leq \tfrac{1}{2}\|X\beta_0 - X\beta\|_{\ell_2}^2 + \langle z, X\hat\beta - X\beta_0\rangle \\ &\quad + 2\lambda_p(\|\beta_0\|_{\ell_1} - \|\hat\beta\|_{\ell_1}). \end{aligned}$$

Put $h = \hat\beta - \beta_0$. We follow the same steps as in Section 3.2 to arrive at

$$\tfrac{1}{2}\|X\hat\beta - X\beta\|_{\ell_2}^2 \leq \tfrac{1}{2}\|X\beta_0 - X\beta\|_{\ell_2}^2 + \langle h_{I_0}, v\rangle - (2 - \sqrt{2})\lambda_p\|h_{I_0^c}\|_{\ell_1},$$

where $v = X^*_{I_0} z - 2\lambda_p \mathrm{sgn}(\beta_{I_0})$. Just as before,

$$\langle h_{I_0}, v \rangle = \langle X^*_{I_0} X h, (X^*_{I_0} X_{I_0})^{-1} v \rangle - \langle h_{I_0^c}, X^*_{I_0^c} X_{I_0} (X^*_{I_0} X_{I_0})^{-1} v \rangle \equiv A_1 - A_2.$$

By assumption, $|A_2| \leq (2 - \sqrt{2})\lambda_p \cdot \|h_{I_0^c}\|_{\ell_1}$. The difference is now in $A_1$, since we can no longer claim that $\|X^* X h\|_{\ell_\infty} \leq (2 + \sqrt{2})\lambda_p$. Decompose $A_1$ as

$$A_1 = \langle X^*_{I_0} X(\hat{\beta} - \beta), (X^*_{I_0} X_{I_0})^{-1} v \rangle + \langle X^*_{I_0} X(\beta - \beta_0), (X^*_{I_0} X_{I_0})^{-1} v \rangle \equiv A_1^0 + A_1^1.$$

Because $\|X^* X(\hat{\beta} - \beta)\|_{\ell_\infty} \leq (2 + \sqrt{2})\lambda_p$, one can use the same argument as before to obtain

$$A_1^0 \leq 2(2 + \sqrt{2})^2 \lambda_p^2 S.$$

We now look at the other term. Since we assume $\|X_{I_0}(X^*_{I_0} X_{I_0})^{-1}\| \leq \sqrt{2}$, we have

$$\begin{aligned}
|A_1^1| &= \langle X(\beta - \beta_0), X_{I_0}(X^*_{I_0} X_{I_0})^{-1} v \rangle \\
&\leq \|X(\beta - \beta_0)\|_{\ell_2} \|X_{I_0}(X^*_{I_0} X_{I_0})^{-1} v\|_{\ell_2} \\
&\leq \sqrt{2} \|X(\beta - \beta_0)\|_{\ell_2} \|v\|_{\ell_2}.
\end{aligned}$$

Using $ab \leq (a^2 + b^2)/2$ and $\|v\|_{\ell_2}^2 \leq (2 + \sqrt{2})^2 \lambda_p^2 S$ gives

$$|A_1^1| \leq \tfrac{\sqrt{2}}{2} \|X(\beta - \beta_0)\|_{\ell_2}^2 + \tfrac{\sqrt{2}}{2}(2 + \sqrt{2})^2 \lambda_p^2 S.$$

To summarize,

$$\langle h_{I_0}, v \rangle \leq \tfrac{\sqrt{2}}{2} \|X(\beta - \beta_0)\|_{\ell_2}^2 + (2 + \tfrac{\sqrt{2}}{2})(2 + \sqrt{2})^2 \lambda_p^2 S + (2 - \sqrt{2})\lambda_p \cdot \|h_{I_0^c}\|_{\ell_1}.$$

It follows that

$$\tfrac{1}{2} \|X\hat{\beta} - X\beta\|_{\ell_2}^2 \leq \tfrac{1 + \sqrt{2}}{2} \|X\beta_0 - X\beta\|_{\ell_2}^2 + (4 + \sqrt{2})(1 + \sqrt{2})^2 \lambda_p^2 S.$$

This concludes the proof.

We close this section by arguing about (1.18) and (1.19). First, it follows from our proof that (1.18) holds. Second, our analysis also shows that the set $\mathcal{A}_{0,S}$ is very large and obeys (1.19).

3.5. *Proof of Theorem 1.3.* Just as with our other claims, we begin by stating a few assumptions that hold with very large probability, and then we show that, under these conditions, the conclusions of the theorem hold. These assumptions are as follows:

  (i) The matrix $X^*_I X_I$ is invertible and obeys $\|(X^*_I X_I)^{-1}\| \leq 2$;
  (ii) $\|X^*_{I^c} X_I (X^*_I X_I)^{-1} \mathrm{sgn}(\beta_I)\|_{\ell_\infty} < \tfrac{1}{4}$;
  (iii) $\|(X^*_I X_I)^{-1} X^*_I z\|_{\ell_\infty} \leq 2\lambda_p$;

(iv) $\|X_{I^c}^*(I - P[I])z\|_{\ell_\infty} \leq \sqrt{2}\lambda_p$;

(v) The matrix–vector product $(X_I^*X_I)^{-1}\mathrm{sgn}(\beta_I)$ obeys

(3.23) $$\|(X_I^*X_I)^{-1}\mathrm{sgn}(\beta_I)\|_{\ell_\infty} \leq 3.$$

We already know that conditions (i) and (ii) hold with large probability [see Section 3.3; the change from $1/2$ to $1/4$ in (ii) is unessential]. As before, we let $E$ be the event $\{\|X_I^*X_I - \mathrm{Id}\| \leq 1/2\}$. For (iii), the idea is the same, and we express $\|(X_I^*X_I)^{-1}X_I^*z\|_{\ell_\infty}$ as $\max_{i \in I} |\langle W_i, z\rangle|$, where $W_i$ is now the $i$th row of $(X_I^*X_I)^{-1}X_I^*$. On $E$, $\max_i \|W_i\| \leq \|(X_I^*X_I)^{-1}X_I^*\| \leq \sqrt{2}$, and the claim now follows from (3.5). Indeed, one can check that conditional on $E$

$$\mathbb{P}(\|(X_I^*X_I)^{-1}X_I^*z\|_{\ell_\infty} > 2\lambda_p) \leq |I| \cdot p^{-2} \cdot (2\pi \log p)^{-1/2}.$$

For (iv), we write $\|X_{I^c}^*(I - P[I])z\|_{\ell_\infty}$ as $\max_{i \in I^c} |\langle W_i, z\rangle|$, where $W_i = (I - P[I])X_i$. We have $\|W_i\| \leq \|X_i\| = 1$ and, conditional on $E$, it follows, from (3.5), that

$$\mathbb{P}(\|X_{I^c}^*(I - P[I])z\|_{\ell_\infty} > \sqrt{2}\lambda_p) \leq |I^c| \cdot p^{-2} \cdot (2\pi \log p)^{-1/2}.$$

The subtle estimate is (v), and it is proven in the next section. There, we show that (3.23) holds with probability at least $1 - 2p^{-2\log 2} - 2|I|\, p^{-2}$. Hence, under the assumptions of Theorem 1.3, (i)–(v) hold with probability at least $1 - 2p^{-1}((2\pi \log p)^{-1/2} + |I|/p) - O(p^{-2\log 2})$.

LEMMA 3.4. *Suppose that the assumptions* (i)–(v) *hold, and assume that* $\min_{i \in I} |\beta_i|$ *obeys the condition of Theorem 1.3. Then, the lasso solution is given by* $\hat{\beta} \equiv \beta + h$ *with*

(3.24) 
$$h_I = (X_I^*X_I)^{-1}[X_I^*z - 2\lambda_p\,\mathrm{sgn}(\beta_I)],$$
$$h_{I^c} = 0.$$

PROOF. The point $\hat{\beta}$ is the unique solution to the lasso functional if

(3.25)
$$X_i^*(y - X\hat{\beta}) = 2\lambda_p\,\mathrm{sgn}(\hat{\beta}_i), \qquad \hat{\beta}_i \neq 0,$$
$$|X_i^*(y - X\hat{\beta})| < 2\lambda_p, \qquad \hat{\beta}_i = 0,$$

and the columns of $X_T$ are linearly independent where $T$ is the support of $\hat{\beta}$. Consider, then, $h$ as in (3.24), and observe that

$$\|h_I\|_{\ell_\infty} \leq \|(X_I^*X_I)^{-1}X_I^*z\|_{\ell_\infty} + 2\lambda_p\|(X_I^*X_I)^{-1}\mathrm{sgn}(\beta_I)\|_{\ell_\infty} \leq 2\lambda_p + 6\lambda_p.$$

It follows that $\|h_I\|_{\ell_\infty} < \min_{i \in I} |\beta_i|$ and, therefore, $\hat{\beta} = \beta + h$ obeys

$$\mathrm{supp}(\hat{\beta}) = \mathrm{supp}(\beta),$$
$$\mathrm{sgn}(\hat{\beta}_I) = \mathrm{sgn}(\beta_I).$$

We now check that $\hat{\beta} = \beta + h$ obeys (3.25). By definition, we have

$$y - X\hat{\beta} = z - Xh = z - X_I(X_I^*X_I)^{-1}[X_I^*z - 2\lambda_p\mathrm{sgn}(\hat{\beta}_I)],$$

since $\beta$ and $\hat{\beta}$ share the same support and the same signs. Clearly,

$$X_I^*(y - X\hat{\beta}) = 2\lambda_p\,\mathrm{sgn}(\hat{\beta}_I),$$

which is the first half of (3.25). For the second half, let $P[I] = X_I(X_I^*X_I)^{-1}X_I^*$ be the orthonormal projection onto the span of $X_I$. Then,

$$\begin{aligned}
\|X_{I^c}^*(y - X\hat{\beta})\|_{\ell_\infty} &= \|X_{I^c}^*(I - P[I])z + 2\lambda_p X_{I^c}^*X_I(X_I^*X_I)^{-1}\mathrm{sgn}(\beta_I)\|_{\ell_\infty} \\
&\leq \|X_{I^c}^*(I - P[I])z\|_{\ell_\infty} + 2\lambda_p\|X_{I^c}^*X_I(X_I^*X_I)^{-1}\mathrm{sgn}(\beta_I)\|_{\ell_\infty} \\
&< \sqrt{2}\lambda_p + \tfrac{1}{2}\lambda_p \\
&< 2\lambda_p.
\end{aligned}$$

Finally, note that $X_T^*X_T$ is indeed invertible, since $T = I$; this is just our invertibility condition. This concludes the proof.  □

Lemma 3.4 proves that $\hat{\beta}$ has the same support as $\beta$ and the same signs as $\beta$, which is of course the content of Theorem 1.3.

3.6. *Proof of (3.23)*. We need to show that $\|(X_I^*X_I)^{-1}\mathrm{sgn}(\beta_I)\|_{\ell_\infty} \leq 3$ with high probability. To begin, we write

$$\begin{aligned}
\|(X_I^*X_I)^{-1}\mathrm{sgn}(\beta_I)\|_{\ell_\infty} &\leq \|\mathrm{sgn}(\beta_I)\|_{\ell_\infty} + \|((X_I^*X_I)^{-1} - \mathrm{Id})\mathrm{sgn}(\beta_I)\|_{\ell_\infty} \\
&\leq 1 + \max_{i \in I}|\langle W_i, \mathrm{sgn}(\beta_I)\rangle|,
\end{aligned}$$

where $W_i$ is the $i$th row of $(X_I^*X_I)^{-1} - \mathrm{Id}$ (or column since this is a symmetric matrix).

LEMMA 3.5. *Let $W_i$ be the $i$th row of $(X_I^*X_I)^{-1} - \mathrm{Id}$. Under the hypotheses of Theorem 1.3, we have*

$$\mathbb{P}\left(\max_{i \in I}\|W_i\| \geq (\log p)^{-1/2}\right) \leq 2p^{-2\log 2}.$$

PROOF. Set $A \equiv \mathrm{Id} - X_I^*X_I$. On the event $E \equiv \{\|\mathrm{Id} - X_I^*X_I\| \leq 1/2\}$ (which holds w.p. at least $1 - p^{-2\log 2}$), we have

$$(X_I^*X_I)^{-1} = I + A + A^2 + \cdots.$$

Therefore, since $W_i = ((X_I^* X_I)^{-1} - \mathrm{Id})e_i$ where $e_i$ is the vector whose $i$th component is 1 and the others 0, $W_i = Ae_i + A^2 e_i + \cdots$ and

$$\|W_i\| \le \|Ae_i\| + \|A\| \|Ae_i\| + \|A^2\| \|Ae_i\| + \cdots$$

$$\le \|Ae_i\| \sum_{k=0}^{\infty} \|A\|^k$$

$$\le \|Ae_i\|/(1 - \|A\|).$$

Hence, on $E$, $\|W_i\| \le 2\|Ae_i\|$.

For each $i \in I$, $Ae_i$ is the $i$th row or column of $\mathrm{Id} - X_I^* X_I$ and for each $j \in I$, its $j$th component is equal to $-\langle X_i, X_j \rangle$ if $j \ne i$, and 0 for $j = i$ since $\|X_i\| = 1$. Thus,

$$\|W_i\|^2 \le 4 \sum_{j \in I : j \ne i} |\langle X_i, X_j \rangle|^2.$$

Now, it follows from Lemma 3.6 that

$$\sum_{j \in I : j \ne i} |\langle X_i, X_j \rangle|^2 \le S\|X\|^2/p + t$$

with probability at least $1 - 2e^{-t^2/[2\mu^2(X)(S\|X\|^2/p + t/3)]}$. Under the assumptions of Theorem 1.3, we have $S\|X\|^2/p \le c_0(\log p)^{-1} \le (8\log p)^{-1}$ provided that $c_0 \le 1/8$. With $t = (8\log p)^{-1}$, this gives

$$(3.26) \qquad \sum_{j \in I : j \ne i} |\langle X_i, X_j \rangle|^2 \le 1/(4\log p)$$

with probability at least $1 - 2e^{-3/[64\mu^2(X)\log p]}$. Now, the assumption about the coherence guarantees that $\mu(X) \le A_0/\log p$ so that (3.26) holds with probability at least $1 - 2e^{-3\log p/[64A_0^2]}$. Hence, by choosing $A_0$ sufficiently small, the lemma follows from the union bound. $\quad\square$

LEMMA 3.6. *Suppose that $I \subset \{1, \ldots, p\}$ is a random subset of predictors with at most $S$ elements. For each $i$, $1 \le i \le p$, we have*

$$(3.27) \qquad \begin{aligned} &\mathbb{P}\bigg(\sum_{j \in I : j \ne i} |\langle X_i, X_j \rangle|^2 > \frac{S}{p}\|X\|^2 + t\bigg) \\ &\le 2\exp\bigg(-\frac{t^2}{2\mu^2(X)(S\|X\|^2/p + t/3)}\bigg). \end{aligned}$$

PROOF. The inequality (3.27) is essentially an application of Bernstein's inequality, which states that, for a sum of uniformly bounded independent random variables with $|Y_k - \mathbb{E}Y_k| < c$,

$$(3.28) \qquad \mathbb{P}\left(\sum_{k=1}^{n}(Y_k - \mathbb{E}Y_k) > t\right) \leq e^{-t^2/(2\sigma^2 + 2ct/3)},$$

where $\sigma^2$ is the sum of the variances, $\sigma^2 \equiv \sum_{k=1}^{n} \mathrm{Var}(Y_k)$. The issue here is that $\sum_{j \in I : j \neq i} |\langle X_i, X_j \rangle|^2$ is not a sum of independent variables and we need to use a kind of Poissonization argument to reduce this to a sum of independent terms.

A set $I'$ of predictors is sampled using a Bernoulli model by first creating the sequence

$$\delta_j = \begin{cases} 1, & \text{w.p.} \ \ S/p, \\ 0, & \text{w.p.} \ \ 1 - S/p, \end{cases}$$

and then setting $I' \equiv \{j \in \{1, \ldots, p\} : \delta_j = 1\}$. The size of the set $I'$ follows a binomial distribution, and $\mathbb{E}|I'| = S$. We make two claims: first, for each $t > 0$, we have

$$(3.29) \qquad \mathbb{P}\left(\sum_{j \in I : j \neq i} |\langle X_i, X_j \rangle|^2 > t\right) \leq 2\mathbb{P}\left(\sum_{j \in I' : j \neq i} |\langle X_i, X_j \rangle|^2 > t\right);$$

second, for each $t > 0$,

$$(3.30) \qquad \begin{aligned} \mathbb{P}&\left(\sum_{j \in I' : j \neq i} |\langle X_i, X_j \rangle|^2 > \frac{S}{p}\|X\|^2 + t\right) \\ &\leq \exp\left(-\frac{t^2}{2\mu^2(X)(S\|X\|^2/p + t/3)}\right). \end{aligned}$$

Clearly, (3.29) and (3.30) give (3.27).

To justify the first claim, observe that

$$\begin{aligned} \mathbb{P}\left(\sum_{j \in I' : j \neq i} |\langle X_i, X_j \rangle|^2 > t\right) &= \sum_{k=0}^{p} \mathbb{P}\left(\sum_{j \in I' : j \neq i} |\langle X_i, X_j \rangle|^2 > t \Big| |I'| = k\right) P(|I'| = k) \\ &\geq \sum_{k=S}^{p} \mathbb{P}\left(\sum_{j \in I' : j \neq i} |\langle X_i, X_j \rangle|^2 > t \Big| |I'| = k\right) P(|I'| = k) \\ &= \sum_{k=S}^{p} \mathbb{P}\left(\sum_{j \in I_k : j \neq i} |\langle X_i, X_j \rangle|^2 > t\right) P(|I'| = k), \end{aligned}$$

where $I_k$ is selected uniformly at random with $|I_k| = k$. We make two observations: (1) since $S$ is an integer, it is the median of $|I'|$ and $P(|I'| \geq S) \geq 1/2$; and (2) $\mathbb{P}(\sum_{j \in I_k : j \neq i} |\langle X_i, X_j \rangle|^2 > t)$ is a nondecreasing function of $k$. To see why this is true, consider that a subset $I_{k+1}$ of size $k + 1$ can be

sampled by first choosing a subset $I_k$ of size $k$ uniformly and then choosing the remaining entry uniformly at random from the complement of $I_k$. It follows that, with $Z_k = \sum_{j \in I_k} |\langle X_i, X_j \rangle|^2 1_{\{i \neq j\}}$, we have that $Z_{k+1}$ and $Z_k + Y_k$, where $Y_k$ is a nonnegative random variable have the same distribution. Hence, $\mathbb{P}(Z_{k+1} \geq t) \geq \mathbb{P}(Z_k \geq t)$. With these two observations in mind, we continue

$$\mathbb{P}\left( \sum_{j \in I': j \neq i} |\langle X_i, X_j \rangle|^2 > t \right) \geq \mathbb{P}\left( \sum_{j \in I: j \neq i} |\langle X_i, X_j \rangle|^2 > t \right) \sum_{k=S}^{p} P(|I'| = k)$$

$$\geq \frac{1}{2} \mathbb{P}\left( \sum_{j \in I: j \neq i} |\langle X_i, X_j \rangle|^2 > t \right),$$

which is the first claim (3.29).

For the second claim (3.30), observe that

$$\sum_{j \in I': j \neq i} |\langle X_i, X_j \rangle|^2 = \sum_{1 \leq j \leq p: j \neq i} \delta_j |\langle X_i, X_j \rangle|^2 \equiv \sum_{1 \leq j \leq p: j \neq i} Y_j.$$

The $Y_j$ are independent and obey:

1. $|Y_j - \mathbb{E}Y_j| \leq \sup_{j \neq i} |\langle X_i, X_j \rangle|^2 \leq \mu^2(X).$
2. The sum of means is bounded by

$$\sum_{1 \leq j \leq p: j \neq i} \mathbb{E}Y_j = \frac{S}{p} \sum_{1 \leq j \leq p: j \neq i} |\langle X_i, X_j \rangle|^2 \leq \frac{S \|X\|^2}{p}.$$

   The last inequality follows from $\sum_{1 \leq j \leq p: j \neq i} |\langle X_i, X_j \rangle|^2 \leq \sum_{1 \leq j \leq p} |\langle X_i, X_j \rangle|^2$ where the right-hand side is equal to $\|X^* X_i\|^2 \leq \|X^*\|^2 \|X_i\|^2 = \|X\|^2$ since the columns are unit-normed.
3. The sum of variances is bounded by

$$\sum_{1 \leq j \leq p: j \neq i} \text{Var}(Y_j) = \frac{S}{p}\left(1 - \frac{S}{p}\right) \sum_{1 \leq j \leq p: j \neq i} |\langle X_i, X_j \rangle|^4 \leq \frac{S \mu^2(X) \|X\|^2}{p}.$$

   The last inequality follows from $\sum_{1 \leq j \leq p: j \neq i} |\langle X_i, X_j \rangle|^4 \leq \mu^2(X) \sum_{1 \leq j \leq p} |\langle X_i, X_j \rangle|^2$, which is less or equal to $\mu^2(X)\|X\|^2$ as before.

The claim (3.30) is now a simple application of Bernstein's inequality (3.27). $\square$

Lemma 3.5 establishes that (3.23) holds with probability at least $1 - 2p^{-2\log 2} - 2|I| \, p^{-2}$. Indeed, on the event $\max_i \|W_i\| \leq (\log p)^{-1/2}$, it follows from Lemma 3.3 that

$$\mathbb{P}\left( \max_{i \in I} |\langle W_i, \text{sgn}(\beta_I) \rangle| \geq 2 \right) \leq 2|I| e^{-2\log p} \leq 2|I| \, p^{-2}.$$

## 4. Discussion.

4.1. *Connection with other works.* In the last few years, there have been many beautiful works that attempt to understand the properties of the lasso and other minimum $\ell_1$ algorithms, such as the Dantzig selector when the number of variables may be larger than the sample size [3, 5, 6, 10, 13, 15, 16, 20, 21, 29, 30]. Some papers focus on the estimation of the parameter $\beta$ and on recovering its support; others focus on estimating $X\beta$. These are quite distinct problems, especially when $p > n$; consider, for instance, the noiseless case.

In [5, 6, 13], it is required that the level of sparsity $S$ be smaller than $1/\mu(X)$. For instance, [5] develops an oracle inequality that requires $S \leq 1/(32\mu(X))$. Even when $\mu(X)$ is minimal [i.e., of size about $1/\sqrt{n}$, as in the case where $X$ is the time-frequency dictionary or about $\sqrt{(2 \log p)/n}$ as for Gaussian matrices and many other kinds of random matrices] one sees that the sparsity level must be considerably smaller than $\sqrt{n}$. When the coherence is of the order of $(\log p)^{-1}$ (as we have allowed in our paper), one would need a sparsity level of order $\log p$. Having a sparsity level substantially smaller than the inverse of the coherence is a common assumption in the modern literature on the subject, although, in some circumstances, a few papers have developed some weaker assumptions. To be a little more specific, [29] reports an asymptotic result in which the lasso recovers the exact support of $\beta$ provided that the strong irrepresentable condition of Section 3.3 holds. The references [20, 28] develop very similar results and use very similar requirements. The recent paper [17] develops similar results but requires either a good initial estimator or a level of coherence on the order of $n^{-1/2}$. In [10, 21] the singular values of $X$ restricted to any subset of size proportional to the sparsity of $\beta$ must be bounded away from zero while [3] introduces an extension of this condition. In nearly all these works, a sufficient condition is that the sparsity be much smaller than the inverse of the coherence.

4.2. *Our contribution.* It follows from the previous discussion that there is a disconnect between the available literature and what practical experience shows. For instance, the lasso is known to work very well empirically when the sparsity far exceeds the inverse of the coherence $1/\mu(X)$ [13], even though the proofs assume that the sparsity is less than a fraction of $1/\mu(X)$. In that paper, the coherence is $1/\sqrt{n}$ so that, as mentioned earlier, results are available only when the sparsity is much smaller than $\sqrt{n}$, which does not explain what series of computer experiments reveal.

Our work bridges this gap. We do so by considering the performance of the lasso one expects in almost all cases but not all. By considering statistical ensembles much as in [9], one shows that, in the above examples, the lasso works provided that the sparsity level is bounded by about $n/\log p$; that is, for generic signals, the sparsity can grow almost linearly with the sample size. We also prove that, under these conditions, the "Irrepresentable Condition" holds with high probability, and we show that, as long as the entries of $\beta$ are not too small, one can recover the exact support of $\beta$ with high probability.

Finally, there does not seem much room for improvement, as all of our conditions appear necessary as well. In Section 2, we have proposed special examples in which the lasso performs poorly. On the one hand, these examples show that, even with highly incoherent matrices, one cannot expect good performance in all cases unless the sparsity level is very small. And on the other hand, one cannot really eliminate our assumption about the coherence, since we have shown that, with coherent matrices, the lasso would fail to work well on generically sparse objects.

One could of course consider other statistical descriptions of sparse $\beta$'s and/or ideal models, and leave this issue open for further research.

**Acknowledgments.** E. Candès would like to thank Chiara Sabatti for fruitful discussions and for offering some insightful comments about an early draft of this paper. E. Candès. also acknowledges inspiring conversations with Terence Tao and Joel Tropp about parts of this paper. We thank the anonymous referees for their constructive comments.

APPLIED AND COMPUTATIONAL MATHEMATICS
CALIFORNIA INSTITUTE OF TECHNOLOGY
300 FIRESTONE, MAIL CODE 217-50
PASADENA, CALIFORNIA 91125
E-MAIL: emmanuel@acm.caltech.edu
        plan@acm.caltech.edu